\documentclass[12pt]{amsart}
\usepackage{amsmath,amsfonts,amssymb,amsthm}
\usepackage{pagesize}
 \newtheorem{theorem}{Theorem}%[section]%[subsection]
 \newtheorem{conjecture}{Conjecture}
  
 \newtheorem{lemma}{Lemma}
 
 \theoremstyle{definition}
 
 \theoremstyle{remark}
 
\def\m1{^{-1}}
\def\ov1{\overline}
\def\gp#1{\langle#1\rangle}

\title{Symmetric Units and Group Identities  in Group Algebras. I}
\author{Victor Bovdi}
\dedicatory{Dedicated to  Professor  L.G.~Kov\'acs on his  70th
birthday}
\address{%
Institute of Mathematics,\newline
\indent University of Debrecen,\newline
\indent 4010 Debrecen, Pf.\ 12,\newline
\indent Hungary}
\address{
Institute of Mathematics and Computer Science,\newline
\indent College of Ny\'\i regyh\'aza\newline
\indent S\'ost\'oi \'ut 31/b,\newline
\indent H-4410 Ny\'\i regyh\'aza, Hungary}

\email{vbovdi@math.klte.hu}


\thanks{The research was supported by OTKA  No.T 037202 and  No.T 038059}
\date{}
\pagesize{paper=a4,width=14cm,height=21.5cm}
\begin{document}
\setcounter{page}{149}
\noindent
{\it Acta Mathematica Academiae Paedagogicae
  Ny\'{\i}regyh\'aziensis}\\
{\bf 22} (2006), 149--159\\
{\tt www.emis.de/journals\\
%\hskip 1.077ex{\small\tt www.emis.de\\\hskip 1.17ex
ISSN 1786-0091}\\[2em]
\allowdisplaybreaks

\begin{abstract}
We describe those   group algebras  over fields of characteristic
different from $2$ whose  units symmetric with respect  to the
classical involution, satisfy some group identity.
\end{abstract}

\maketitle
%%%%%%%%%%%%%%%%%%%%%%%%%%%%%%%%%%%%%
\section{Introduction}

Let  $U(A)$ be the group of units of an algebra  $A$ with
involution $*$ over the field $F$ and let $S_*(A)=\{u\in U(A)\mid
u=u^*\}$ be the set of symmetric units of $A$.

Algebras with involution have been actively investigated. In these
algebras there are many symmetric elements, for example: $x+x^*$
and $xx^*$  for any $x\in A$. This raises  natural questions about
the properties of the symmetric elements and symmetric units. In
\cite{10} Ch.~Lanski  began to study the properties of the
symmetric units in prime algebras with involution, in particular
when the symmetric units commute. Using the results and methods of
\cite{3}, in \cite{4} we classified the cases when the symmetric
units commute in  modular group algebras of $p$-groups. The
solution of this question for integral group rings and for some
modular group rings of arbitrary groups was obtained in \cite{5,
6}.

Several results on the group of units $U(R)$ show  that if $U(R)$
satisfies a certain group theoretical condition (for example, it
is  nilpotent or solvable), then $R$'s properties  are  restricted
and a polynomial identity on  $R$ holds. This suggests that there
may be some general underlying relationship between group
identities and polynomial  identities. In this topic Brian Hartley
made the following:
\begin{conjecture} Let $FG$ be a group algebra
of a torsion group $G$ over the field $F$.
 If $U(FG)$ satisfies a group identity, then $FG$ satisfies a polynomial
 identity.
\end{conjecture}

The theory of $PI$-algebras has been established  for a long time.
On the contrary, the  study of algebras  with units satisfying a
group identity has emerged only recently \cite{11, 12}. Our goal
here  is to show that with a  few  extra assumptions, these
algebras are actually $PI$-algebras. In fact, these classes of
algebras are quite special, because if the group of units is too
small in a algebra, a group identity condition can not limit the
structure of the whole algebra. In view of Hartley's conjecture,
as a natural generalization the works \cite{4, 5, 6, 10} it is a
natural question when does the symmetric units satisfy a group
identity in group algebra. Note that the structure theorem of the
algebras with involution whose symmetric elements satisfy a
polynomial identity was obtained earlier  by S.A.~Amitsur in
\cite{1}. A.~Giambruno, S.K.~Sehgal and A.~Valenti in  \cite{8}
obtained the following  result for group algebras of torsion
groups:

\begin{theorem} Let $FG$ be a group algeb\-ra of a torsion group
$G$ over an infinite field $F$ of characteristic $p>2$ and assume
that the involution $*$ on $FG$ is canonical. The symmetric units
$S_*(FG)$ satisfy a group identity   if and only if $G$ has a
normal subgroup $A$ of finite index, the commutator subgroup $A'$
is a finite $p$-group and one of the following conditions holds:
\item{(i)} $G$ has no  quaternion subgroup  of order $8$ and $G'$
has  of bounded exponent $p^k$ for some $k$. \item{(ii)} $G$ has
of bounded exponent $4p^s$ for some $s\geq 0$, the $p$-Sylow
subgroup of $G$ is normal  and $G/P$ is a Hamiltonian
$2$-subgroup.
\end{theorem}

\bigskip

In the present paper we extend the result of A.~Giambruno,
S.K.~Sehgal and A.~Valenti. For non-torsion groups $G$ we describe
the  group algebras $FG$ over the field $F$ of characteristic
different from $2$ whose symmetric units
$$
\textstyle S_*(FG)=\{u=\sum_{g\in G}\alpha_gg\in
U(FG)\quad\mid\quad u= u^*=\sum_{g\in G}\alpha_gg^{-1}\}
$$
satisfy  a group identity. The present result was announced at the
International Workshop "Polynomial Identities in Algebras, 2002,
Memorial University of Newfoundland.

\section{Main results}
In the  sequel of this paper  $ \mathfrak{d}(\omega)$ denotes a
positive integer, which depends on the group identity $\omega$ and
it is defined  in the next section. Our results are the following:

\begin{theorem} Let $G$ be a non-torsion nonabelian group  and
 $char(F)=p\not=2$ and assume that the symmetric units of $FG$
satisfy some  group identity $\omega=1$. Assume that $|F|>
\mathfrak{d}(\omega)$, where $ \mathfrak{d}(\omega)$ is an integer
which depends only on the word $\omega$. Let $P$ be a $p$-Sylow
subgroup of $G$ and let $t(G)$ be the torsion part  of $G$.
\begin{itemize}
\item[(I)] If $p>2$ then $P$ and $t(G)$ are normal subgroups of
$G$ such that:
\begin{itemize}
\item[(a)] $B=t(G)/P$ is an abelian $p'$-subgroup and its
subgroups are normal in $G$; \item[(b)] if $B$ is noncentral in
$G/P$ then the algebraic closure $L$ of the prime subfield $F_p$
in $F$ is finite and for all $g\in G/P$ and for any $a\in B$ there
exists an $r\in\mathbb{N}$ such that $a^g=a^{p^r}$ and $|L:F_p|$
is a divisor of $r$; \item[(c)] the $p$-Sylow subgroup $P$  is a
finite group; \item[(d)] the $p$-Sylow subgroup $P$ is infinite
and $G$ has a subgroup $A$ of finite index, such that $A'$ is  a
finite $p$-group and the commutator subgroup $H'$ of $H=AP$ is a
bounded $p$-group. Moreover, if $P$ is unbounded, then  $G'$ is a
bounded $p$-group;
\end{itemize}
\item[(II)]  If $char(F)=0$  then $t(G)$ is a subgroup, every
subgroup of $t(G)$ is  normal in $G$ and one of the following
conditions holds:
\begin{itemize}
\item[(a)] $t(G)$ is abelian and each idempotent of $Ft(G)$  is
central in $FG$; \item[(b)]  $t(G)$ is a Hamiltonian $2$-group,
and each symmetric idempotent of $Ft(G)$ is  central in $FG$.
\end{itemize}
\end{itemize}
\end{theorem}

\section{Notation, preliminary results and the proof}

Let $FG$ be the group algebra of $G$ over $F$. We introduce the
following notation:
\begin{itemize}
\item[$\bullet$] $(g,h)=g\m1g^h=g\m1h\m1gh$ for all $g,h\in G$;
\item[$\bullet$] $|g|$ and $C_G(g)$ are the order and the
centralizer of  $g\in G$, respectively; \item[$\bullet$] $G'$,
$Syl_p(G)$ are the commutator subgroup and the Sylow $p$-subgroup
of $G$; \item[$\bullet$] $t(G)$ is the set of elements of finite
order in  $G$; \item[$\bullet$] $\Delta(G)= \{ g\in G \mid \quad
[G:C_G(g)]<\infty \}$ is the $FC$-radical of $G$; \item[$\bullet$]
$\Delta^p(G)= \gp{ g\in \Delta(G)  \mid \quad g \text{ has order
of a power of p} }$; \item[$\bullet$] $T_{l}(G/H)$ is a left
transversal of the  subgroup $H$ in  $G$; \item[$\bullet$] $
\mathfrak{N}(FG)$ is the sum of all nilpotent ideals of the group
algebra $FG$; \item[$\bullet$] $A(FG)$ is the augmentation ideal
of the group algebra $FG$.
\end{itemize}

\bigskip
Let  $A$ be an algebra over a field $F$, let $F_0$ be the ring of
integers of the field $F$, and suppose  that $U(A)$ satisfies a
group identity $\omega=1$. Then, as it was proved in Lemma $3.1$
of  \cite{11}, there exists a polynomial $f(x)$ over $F_0$ of
degree $ \mathfrak{d}(\omega)$ which is determined by the word
$\omega$. In several papers (see for example \cite{8}) the authors
assumed that the field $F$ is infinite so  they could apply the
``Vandermonde determinant argument''. We shall use some  lemmas from
\cite{8}, which are easy to prove using the method of the paper
\cite{11} even without the assumption  that the field  $F$ is
infinite.

In our proof we will use  the following facts:

\begin{lemma}\label{l1}(\cite{1}) Let $A$ be an  algebra with
involution over  $F$ of  $char(F)\ne 2$, such that the set of
symmetric units of $A$ satisfy a group identity $\omega=1$. If $I$
is a stable nil ideal of $A$ then the  symmetric units  of $A/I$
satisfy a group identity.
\end{lemma}

\begin{lemma}\label{l2}(see \cite{8}) Let $A$ be an  algebra over the
field $F$ of   characteristic $p\ne 2$, such that the set of
symmetric units  of $A$ satisfy a group identity $\omega=1$   and
$|F|> \mathfrak{d}(\omega)$, where $ \mathfrak{d}(\omega)$ is an
integer which depends only on the word $\omega$. Then:
\begin{itemize}
\item[(i)]  if $A$ is semiprime, then $asa=0$ for every nilpotent
element  $s\in S_*(A)$ and square-zero $a\in S_*(A)$; \item[(ii)]
if $a\in A$ is  square-zero, then $(aa^*)^m=0$, for some $m\in
\mathbb{N}$; \item[(iii)] if $A$ is semiprime  and $u,v\in A$ such
that $uv=0$, then  $usv=0$ for any square-zero symmetric element
$s$; \item[(iv)] if the subring $L$  of $A$ is nil, then $L$
satisfy a polynomial identity; \item[(v)] each symmetric
idempotent is central; \item[(vi)]  if $A$ is artinian, then $A$
is isomorphic to a direct sum of division algebras and $2\times 2$
matrices algebras over a field with symplectic involution. Each
nilpotent element of $A$ has index at most $2$; \item[(vii)]  if
$A=FG$  is  the group algebra of the group $G=Q_8\times \gp{c}$,
where $Q_8$ is the quaternion group of order $8$, then the order
of   the cyclic subgroup $\gp{c}$ is finite.
\end{itemize}
\end{lemma}

\begin{lemma}\label{l3}(see \cite{8}) Let $A$ be a normal abelian
subgroup of $G$ of finite index such that   $G=A\cdot H$, where
$H$ is a finite group. Let  $char(F)=p$ and  assume that the set
of symmetric units of $FG$ satisfy a group identity $\omega=1$. If
$|F|> \mathfrak{d}(\omega)$, where $ \mathfrak{d}(\omega)$ is an
integer which depends only on the word $\omega$, then $G'$ has
bounded exponent $p^m$, where $m$ depends only on $ \mathfrak{d}$.
\end{lemma}

Now we are ready to prove  the following

\begin{lemma}\label{l4} Let  $char(F)=p>2$ and  let the set of
symmetric units  of $FG$ satisfy a group identity $\omega=1$.
Assume that $|F|> \mathfrak{d}(\omega)$, where $
\mathfrak{d}(\omega)$ is an integer which depends only on the word
$\omega$. Then the $p$-Sylow subgroup $P$ of $\Delta (G)$ is
normal in  $G$ and  the set of symmetric units of $F[G/P]$ satisfy
a group identity.
\end{lemma}

\begin{proof} Let $H$ be a finite subgroup  of $\Delta(G)$ and let
$J=J(F_pH)$ be the radical of the finite group algebra $F_pH$ over
the prime subfield $F_p$. According to  Lemma {\sl{2(vi)}}, the
factor algebra $F_pH/J$ is isomorphic to a direct sum of fields
and $2\times 2$ matrices algebras over a finite field with
symplectic involution and  a nilpotent element $\bar u=u+J\in
F_pH/J$ has index at most $2$. Moreover, from  this decomposition
follows that $\bar u\bar u^*$ is central.  By Lemma 2{(\sl{ii})}
the element $\bar u\bar u^*$ is nilpotent and central in the
semiprime algebra $F_pH/J$.  Therefore $\bar u\bar u^*=0$ and
$uu^*\in J(FH)$.

Let $h\in H$ with   $|h|=p^t$. Then  $u=h-1$ is nilpotent and
$$
uu^*=(h-1)(h^{-1}-1)\in J(FH).
$$
It follows that $huu^*=-(h-1)^2\in J(FH)$. Using Passman's result
(see Lemma 5 in \cite{8}, p.453) we obtain that $h-1\in J(FH)$ for
all $h\in H$ and $H\cap (1+J)$ is a normal $p$-subgroup of $H$,
which coincides with the $p$-Sylow subgroup of $H$. Thus  the
$p$-Sylow subgroup $P$ of $\Delta (G)$ is normal in  $G$, so the
proof is complete.\qquad
\end{proof}

\bigskip

\begin{lemma}\label{l5} Let $FG$ be a semiprime group algebra over
the field $F$ with\linebreak $char(F)>2$  such that the set of symmetric
units  of $FG$ satisfy a group identity $\omega=1$. Suppose that
$|F|> \mathfrak{d}(\omega)$, where $ \mathfrak{d}(\omega)$ is an
integer which depends only on the word $\omega$. Then one of the
following conditions holds:
\begin{itemize}
\item[(i)] $t(G)$ is abelian and each idempotent of $Ft(G)$  is
central in $FG$. \item[(ii)]  $t(G)$ is a Hamiltonian $2$-group
and each symmetric idempotent of $Ft(G)$ is central in $FG$.
\end{itemize}
\end{lemma}

\begin{proof} {(i)}  Let $a\in t(G)$, such that $(|a|,p)=1$. Then,
by Lemma $2${\sl{(v)}},  the symmetric idempotent\quad
$e=\frac{1}{n}(1+a+\cdots +a^{|a|-1})$\quad  is central in $FG$,
so $\gp{a}$ is normal in $G$. Now let $p>2$ and let $a\in t(G)$ be
of order $p$. If $N_G(\gp{a})=G$ then $\overline{\gp{a}}$ is a
central nilpotent element of the semiprime algebra $FG$, a
contradiction.

Let us  prove that each torsion element  belongs to $N_G(\gp{a})$.
Pick  $h\not\in N_G(\gp{a})$ such that $|h|=p^t$. The elements
$(h-1)(h^{-1}-1)$ and  $\overline{\gp{a}}$ are symmetric and \quad
$(2-h-h^{-1})^{p^t}=(\overline{\gp{a}})^2=0$.\quad  By Lemma
$2${(\sl{i})} we get
$\overline{\gp{a}}(2-h-h^{-1})\overline{\gp{a}}=0$ and
\begin{equation}\label{e1}
\overline{\gp{a}}h\overline{\gp{a}}+
\overline{\gp{a}}h^{-1}\overline{\gp{a}}=0.
\end{equation}
An   element of $Supp(\overline{\gp{a}h}\overline{\gp{a}})$ can be
written  as $a^{i}ha^j$, where $0\leq i,j\leq p-1$. If all the
elements in    $Supp(\overline{\gp{a}}h\overline{\gp{a}})$ and in
$Supp(\overline{\gp{a}}h^{-1}\overline{\gp{a}})$ are distinct,
then on the left-hand side of (\ref{e1}) each element appears at
most two times, but this leads to a contradiction if
$char(F)\not=2$. Therefore, in the subset
$Supp(\overline{\gp{a}}h\overline{\gp{a}})$  not all elements are
different, whence there exist $i,j,k,l$ such that
$a^iha^j=a^kha^l$ and either $i\ne k$ or $j\ne l$. If, for
example,  $i>k$, then $h^{-1}a^{i-k}h=a^{l-j}$ and $h\in
N_G(\gp{a})$.

Now, let $h\not\in N_G(\gp{a})$  be a $p'$-element. As we have
seen before, $\gp{h}$ is normal in  $G$, so   $\gp{a,h}$ is a
finite subgroup.  By Lemma $4$  the $p$-Sylow subgroup $P$ of
$\gp{a,h}$ is normal in $\gp{a,h}$ and $(a,h)\in P\cap
\gp{h}=\gp{1}$, a contradiction.

Therefore, each  element  of finite  order  belongs to
$N_G(\gp{a})$. Moreover,  the  elements of order $p$ in $G$ form
an elementary abelian normal  $p$-subgroup $E$ of $G$.

Finally,  if  $h\not\in N_G(\gp{a})$, then $h$ has  infinite order
and $h$ acts on $E$. The subgroups   $\gp{a^{h}}$ and  $\gp{a}$
are different  and we can choose a subgroup $\gp{b}\subset E$,
which  differs from $\gp{a}$. Clearly,
$\overline{\gp{a}}(h+h^{-1})\overline{\gp{a}}$ and
$\overline{\gp{b}}$ are square-zero symmetric elements and
according to Lemma $2${(\sl{i})},
\begin{equation}\label{e2}
\overline{\gp{b}} \overline{\gp{a}}(h+h^{-1})\overline{\gp{a}}
\overline{\gp{b}}=0.
\end{equation}

Since  $hE$ and $h^{-1}E$ are different cosets,  from (\ref{e2})
follows that
\begin{equation}\label{e3}
\overline{\gp{b}}
\overline{\gp{a}}h\overline{\gp{a}}\overline{\gp{b}}=0.
\end{equation}

The subgroup $H=\gp{a,b}\subset E$ has order $p^2$ and by
(\ref{e3}) we have  ${\overline H}h_1{\overline H}h_2=0$ for all
$h_1,h_2$. Since elements of finite order belong to $N_G(H)$, we
get $({\overline H} FG)^2=0$, which is impossible by the
semiprimeness of $FG$. Thus $G$ have no  $p$-elements and  all
finite cyclic subgroups of $G$ are normal in $G$. Applying Lemmas
$6$ and $7$ from \cite{8} and the  fact that $G$ have no
$p$-elements ($p\not=2$), we obtain that $t(G)$ is either an
abelian group or a Hamiltonian $2$-group.

Let $t(G)$ be an abelian group and let $e\in Ft(G)$ be a
noncentral idempotent in $FG$. Set $H=\gp{Supp(e)}$. Since every
subgroup of $t(G)$ is normal in $G$, the subgroup $H$ is also
normal in $G$ and $FH$ has a primitive idempotent $f$, which  does
not commute with some $g\in G$ of infinite order. Then $g\m1
fg\not=f$ is also a primitive idempotent of $FH$ and
$(g^{-1}fg)f=0$,  i.e. $(fg)^2=(gf)^2=0$.

Let $g\m1 fg=\ov1{f}\not=f^*$. By Lemma $2$(v) we have\quad
$f\not=f^*$, so\quad   $g\m1f+f^*g$\quad  is a square-zero
symmetric element and by  Lemma $2${\sl(iii)}, we get that
$$
fg(g\m1f+f^*g)fg=0.
$$
It follows that\quad  $f+g(\ov1{f}f^*)gf=f=0$,\quad  a
contradiction. Therefore,\quad  $g\m1fg=f^*$, so\quad
$(f^*)^*=(g\m1fg)^*=g\m1f^*g=f$. Furthermore,
$g^{-2}fg^2=g\m1f^*g=f$\quad and\quad  $f^*g^2=g^2f^*$. Since
$f^*g^2=g^2f^*$,\quad  $(gf^*)^2=0$ and $gf+f^*g\m1$ is
square-zero symmetric element, by  Lemma $2${\sl (iii)}  we obtain
that
$$
gf^*(gf+f^*g\m1)gf^*= gf^*g^2(g\m1fg)f^*+gf^*=gf^*g^2f^*+gf^*=0.
$$
Thus  $(g^2+1)f^*=0$, which is impossible, since $g^2H$ and $H$
are different cosets.
\end{proof}

\begin{lemma}\label{l6} Let $F$ be a field of characteristic $p$, and
suppose that $G$ contains a normal locally finite $p$-subgroup $P$
such that the centralizer of each element of $P$ in every finitely
generated subgroup of $G$ is of finite index. Then  $
\mathfrak{I}(P)$ is a locally nilpotent ideal.
\end{lemma}

\begin{proof} Clearly,  $\{\; u(h-1) \;\mid\; u \in T_{l}(G/P),\; 1
\neq h \in P\;\} $ is   an $F$-basis for the ideal ${ \mathfrak{
I}}(P)$. Let us show that the subalgebra $W=\gp{u_1(h_1-1), \;
\ldots, \; u_s(h_s-1)}_F$ is nilpotent. According to our
assumption, the centralizers of $h_1,\ldots, h_s$ in the subgroup
$H=\langle u_1,\ldots, u_s,h_1,\ldots,h_s\rangle$ have  finite
index. Since $P$ is normal, its  subgroup $ L=\langle h_1^u,
h_2^u,\ldots, h_s^u\mid u\in H\rangle $ is a finitely generated
FC-group and by a Theorem of B.H.~Neumann (\cite{1}, Theorem 4,
p.19) $L$ is a finite $p$-group. Thus the augmentation ideal
$A(FL)$ is nilpotent with index, say, $t$. Furthermore,
$A(FL)=u^{-1}A(FL)u$ for any $u\in H$ and this implies that
$\left(A(FL)\cdot FH\right)^n=A^n(FL)\cdot FH$  for any $n>0$, so
$W^t\subseteq A^t(FL)\cdot FH=0$, because $W\subseteq A(FL)\cdot
FH$. Therefore $W$ is a nilpotent subalgebra and ${
\mathfrak{I}}(P)$ is a locally nilpotent ideal.\qquad

\end{proof}

\begin{lemma}\label{l7} Let    $G$ be a group with a nontrivial
$p$-Sylow subgroup $P$ and let $char(F)=p>2$. If the set of
symmetric units  of $FG$ satisfy a group identity $\omega=1$   and
$|F|> \mathfrak{d}(\omega)$, where $ \mathfrak{d}(\omega)$ is an
integer which depends only on the word $\omega$, then  $P$ is
normal in $G$ and the ideal ${ \mathfrak{I}}(P)$ is nil.
\end{lemma}

\begin{proof}  Let $P$ be a maximal normal $p$-subgroup of $G$ such
that the ideal ${ \mathfrak{I}}(P)$ is  nil. By Lemma $1$  the set
of symmetric units  of $F[G/P]$ satisfy a group identity. If
$F[G/P]$ is semiprime, then by (i) of the  Theorem  the group
$G/P$ has no $p$-elements and $P$ coincides with the $p$-Sylow
subgroup of $G$. Now, suppose that $F[G/P]$ is not semiprime.
According to Theorem 4.2.13 (\cite{13}, p.131)  the group
$\Delta(G/P)$ has a nontrivial $p$-Sylow subgroup $P_1/P$, which
is normal in $G/P$ by Lemma $4$. Clearly, $P_1/P$ is an
$FC$-subgroup of $G/P$, so by Lemma $6$ the ideal ${
\mathfrak{I}}(P_1/P)$ is nil.

Since ${ \mathfrak{I}}(P_1/P)\cong { \mathfrak{I}}(P_1)/{
\mathfrak{I}}(P)$ and $P_1$ is normal in $G$, the ideal  ${
\mathfrak{I}}(P_1)$ is nil and $P\subset P_1$, a contradiction.
Thus $P=Syl_p(G)$ and the proof is done. \qquad
\end{proof}

\begin{lemma}\label{l8} Let $R$ be an algebra with involution $*$ over a
field $F$ of characteristic $p>2$   and assume that $S_*(R)$
satisfies a group identity and $|F|> \mathfrak{d}(\omega)$. If
some nil subring $L$ of $R$ is $*$-stable, then $L$ satisfies  a
non-matrix polynomial identity.
\end{lemma}

\begin{proof}  Let  $A=F\gp{X }[[t]]$ be  the ring of   power
series over the  polynomial ring $F\gp{X}$ with noncommuting
indeterminably $X=\{ x_1,x_2 \}$. By a  result  of Magnus, the
elements $1+x_1t,1+x_2t$ are units in $A$ and $\gp{1+x_1t,
1+x_2t}$ is a free group.

Assume that $S_*(R)$  satisfies the  group identity $w$, where $w$
is a reduced word in  $2$ variables. Then $w(1+x_1t, 1+x_2t)\ne 1$
according to result of Magnus and it is well-known that
$(1+x_it)^{-1}=1-x_it+x_i^2t^2-\cdots$. If we substitute
$(1+x_it)^{-1}$ in the expression $w(1+x_1t,1+x_2t)-1$, then it
can be  expanded  as
\begin{equation}\label{e4}
\textstyle \sum_{i\geq s} g_i(x_1,x_2)t^i,
\end{equation}
where $g_i(x_1,x_2)\in F\gp{X}$ is a  homogeneous polynomial of
degree $i$.  Obviously  there exists a smallest integer $s\geq 1$
such that   $g_s(x_1,x_2)\ne 0$.

Let $L$ be a $*$-stable nil subring and  let $S(L)$ be the set of
the  symmetric  elements of $L$. Take now $r_1, r_2\in S(L)$ and
let $\lambda\in F$. Obviously, $r_1, r_2$ are nilpotent elements,
so each  $1+\lambda r_i$ is a symmetric unit in $R$ and
$$
(1+r_i\lambda)^{-1}=1-r_i\lambda+r_i^2\lambda^2+\cdots
+(-1)^{t-1}r_i^{t-1}\lambda^{t-1}
$$
for a suitable $t$. By evaluating $w$ on these elements,
(\ref{e4}) gives us a finite sum $ \textstyle \sum_{i\geq s}^l
g_i(r_1,r_2)\lambda^i=0$ for some $l$. Since $|F|>
\mathfrak{d}(\omega)$, we can apply the  Vandermonde determinant
argument to obtain $g_i(r_1,r_2)=0$ for all $i$. Therefore
$g_s(x_1,x_2)$ is a $*$-polynomial identity on $S(L)$. Finally, by
\cite{1} it follows that $S(L)$ satisfies an ordinary polynomial
identity.

Suppose  that   the  homogeneous polynomial  $g(x_1,x_2)$ vanishes
on the matrix algebra  $M_2(K)$ over a commutative ring $K$. Then
$$
g(x_1,x_2)
=h(x_1,x_2)+g_{11}(x_1,x_2)+g_{12}(x_1,x_2)+g_{21}(x_1,x_2)+
g_{22}(x_1,x_2),
$$
where $h(x_1,x_2)$ consists of all monomials which contain $x_1^2$
or $x_2^2$  while the  $g_{ij}(x_1,x_2)$ contain all the remaining
monomials beginning  with $x_i$ and  ending  with $x_j$ for
$i,j\in \{1,2\}$. If $a$ and $b$ are two square-zero matrices,
then $h(a,b)=0$, because each term of $h$ has $a^2$ or $b^2$ as a
factor, so  we conclude that  $ag_{21}(a,b)b=0$. Clearly
$x_1g_{21}(x_1,x_2)x_2$ is some polynomial $f(x_1x_2)$. Then
$f(ab\lambda)=0$  for each  $\lambda \in F$ and, by the
Vandermonde determinant argument, we get  $(ab)^d=0$ for some $d$.
Take, for instance, the matrix units  $a=e_{12}$ and $b=e_{21}$,
then we obtain a contradiction. \qquad
\end{proof}

\begin{lemma}\label{l9} Let $R$ be an  algebra  over a  field $F$ of
positive characteristic  $p$ satisfying  a non-matrix polynomial
identity. Then  $R$ satisfies also  a  polynomial  identity of the
form  $([x,y]z)^{p^l}$ and  $[x,y]^{p^l}$
\end{lemma}
\begin{proof}\quad Let $g(x_1,x_2,\ldots, x_n)$ be a non-matrix
polynomial identity in $R$. The variety $W$ determined by the
polynomial identity $g(x_1,x_2,\ldots, x_n)$ contains  a
relatively free algebra $K$  of rank 3. Of course, $K$ is a
finitely  generated  $PI$-algeb\-ra,  and the result of  Braun and
Razmyslov (Theorem  6.3.39, \cite{14}) states that the radical
$J(K)$ of $K$ is nilpotent. Writing $K/J(K)$ as a subdirect sum of
primitive rings $\{L_i\}$, we get   that  every primitive ring
$L_i$ satisfies the non-matrix  polynomial identity
$g(x_1,x_2,\ldots, x_n)$, as a  homomorphic image of $K$. By
Theorem 2.1.4 of \cite{9}, $L_i$ is either isomorphic to the
matrix ring $M_m(D)$ over a division ring  $D$,  or for any $m$
the matrix ring $M_m(D)$ is an epimorphic image of some subring of
$L_i$.

Thus   $M_m(D)$ satisfies a non-matrix  polynomial identity  $g$,
which is  possible only if $L_i$ is a commutative ring.
Consequently,  $K/J(K)$ is a commutative algebra, so $K$ satisfies
a polynomial identity  of the form $([x,y]z)^{p^l}$ such that
$J(K)^{p^l}=0$.  Since  $R$ belongs to the variety $W$, the
algebra  $R$ also satisfies a polynomial  identity
$([x,y]z)^{p^l}$. \qquad
\end{proof}

\begin{lemma}\label{l10} Let $FG$ be a non semiprime group algebra over
the field $F$ with $char(F)>2$, such that the set of symmetric
units  of $FG$ satisfy a group identity $\omega=1$   and $|F|>
\mathfrak{d}(\omega)$, where $ \mathfrak{d}(\omega)$ is an integer
which depends only on the word $\omega$. If $ \mathfrak{N}(FG)$ is
not nilpotent   then $FG$ is a $PI$-algebra, where $
\mathfrak{N}(FG)$ is the sum  of all nilpotent ideals of $FG$.
\end{lemma}

\begin{proof}   Clearly the  non nilpotent ideal
$ \mathfrak{N}= \mathfrak{N}(FG)$ is invariant under the
involution $*$ and by  Lemma $2${(\sl{iv})} the ring $
\mathfrak{N}$ satisfies a polynomial identity $f(x_1,\ldots,x_n)$.
Moreover, by Lemma $2.8$ of \cite{12} the algebra $FG$ satisfies a
nondegenerate multilinear generalized polynomial identity and
hence, by Theorem 5.3.15 (\cite{13}, p.202),
$|G:\Delta(G)|<\infty$ and $\Delta(G)'$ is finite.

Set $P=Syl_p(G)$ and $P_1=Syl_p(\Delta(G)')$. By Lemma $4$,
$P\cap\Delta(G)'=P_1\vartriangleleft G$ and $P_1$ is a finite
$p$-group. Thus $ \mathfrak{I}(P_1)$ is a nilpotent ideal and by
(i) of the Theorem, the set of symmetric units  of
$F[\Delta(G)'/P_1]$ satisfy a group identity, so $\Delta(G)'/P_1$
is either an abelian $p'$-group or a Hamiltonian $2$-group.

If $P_1=\Delta(G)'$,  then by Theorem 5.3.9 (\cite{13}, p.197) the
algebra $FG$ is a $PI$-algebra. If $P_1\subsetneq\Delta(G)'$ then
we can suppose that $G$ is  a group such that
$Syl_p(\Delta(G)')=1$ and $\Delta(G)'$    is either an abelian
$p'$-group or a Hamiltonian $2$-group.

Set $P_2=Syl_p(\Delta(G))$. Clearly, $P_2=P\cap \Delta(G)$ is
normal in $\Delta(G)$. Since $[P:P_2]<\infty$ and $P$ is an
infinite group, the group  $P_2$ is infinite,  too. If $a\in P_2$,
$b\in \Delta(G)$, then $(a,b)\in P_2\cap\Delta(G)'=1$, so
$(a,b)=1$ and $P_2$ is a central subgroup in $\Delta(G)$.

Let us prove that   $F\Delta(G)$ is a $PI$-algebra. If $\Delta(G)$
is a torsion  group, then by \cite{8} the statement is trivial.

Since $ \mathfrak{N}(F\Delta(G))\subseteq  \mathfrak{N}(FG)$, the
ideal $ \mathfrak{N}(F\Delta(G))$ also satisfies the same
polynomial identity $f(x_1,\ldots,x_n)$. By the standard
multilinearization process, we may assume that $f(x_1,\ldots,x_n)$
is multilinear.

Assume  that $P_2$ has  bounded exponent.  Then  the maximal
elementary abelian $p$-subgroup  $E$ of $P_2$ is infinite.
Let\quad $f(a_1,\ldots,a_n)=\sum_i \alpha_i y_i$,\quad  where
$a_1,\ldots, a_n\in F\Delta(G)$,\quad   $y_1,\ldots,y_n\in
T_l(\Delta(G)/E)$\quad  and $\alpha_i\in FE$. Then there exists a
finite subgroup $B$ such that $\alpha_i\in FB$ and $E=B\times
\prod_j\gp{c_j}$. Since $(c_k-1)a_k\in  \mathfrak{N}(F\Delta(G))$
and $P_2$ is central, we conclude  that
$$
f((c_1-1)a_1,\ldots,(c_n-1)a_n) =(c_1-1)\cdots
(c_n-1)f(a_1,\ldots,a_n)=0.
$$
It  follows that $f(a_1,\ldots,a_n)=0$, because $B\cap
\prod_j\gp{c_j}=\gp{1}$.

Now let $P_2$ be  of unbounded exponent and $c\in P_2$. Then
$(c-1)a_k\in  \mathfrak{N}(F\Delta(G))$ and also
$$
f((c-1)a_1,\ldots,(c-1)a_n)=(c-1)^nf(a_1,\ldots,a_n)=0
$$
for all $c\in P_2$. Then $f(a_1,\ldots,a_n)$ belongs to the
annihilator of  the augmentation ideal $A(FP_2^{p^t})$, where
$n\leq p^t$. Since   $P_2^{p^t}$ is infinite, we have
\[
Ann_l(A(FP_2^{p^t}))=0.
\]
It follows that $f(a_1,\ldots,a_n)=0$,
so  $f(x_1,\ldots,x_n)$ is a polynomial identity for $F\Delta(G)$.
Since  $F\Delta(G)$ is a $PI$-algebra and $[G:\Delta(G)]<\infty$,
the algebra   $FG$ is $PI$, too.\qquad  \end{proof}

\begin{proof}[Proof of the theorem] Let $FG$ be a group  algebra of
a non-torsion group $G$ over a field of positive characteristic
$p$. By Lemma $7$ the $p$-Sylow subgroup $P$ is normal in $G$ and
$F[G/P]\cong FG/{ \mathfrak{I}}(P)$, so the symmetric units  of
semiprime algebra $F[G/P]$ satisfy a group identity. By Lemma $5$
$B=t(G/P)$ is a subgroup of $G/P$ and $B$ is either an abelian
$p'$-group or a Hamiltonian $2$-group. If $B$ is a Hamiltonian
$2$-group, then $Q_8$ is a subgroup of $B$. Choose an element
$c\in G/P$  of infinite order. Since every subgroup of $t(G)/P$ is
normal in $G/P$ and \quad $|Aut(Q_8)|< \infty$,\quad there exists
a $t\in \mathbb{N}$ such that $c^t\in C_{G/P}(Q_8)$ and $Q_8\times
\gp{c^t}\subseteq G/P$. Then Lemma $2$(vii) asserts that $c$ has
finite order, a contradiction. So $B$ is an abelian $p'$-group and
by Lemma 5 every idempotent of $FB$ is central in $F[G/P]$.
Moreover, if $B$ is noncentral, then according to \cite{7} the
group $B$ satisfy (i.b) of our Theorem.

Now, let  $P$ be infinite. By Corollary 8.1.14 (\cite{13}, p.312)
the ideal $ \mathfrak{N}(FG)$ is  non-nilpotent, so by  Lemma
$10$, the algebra $FG$ is a $PI$-algebra, i.e. $G$ has a subgroup
$A$ with finite index such that $A^\prime$ is a finite $p$-group.
According to  Lemma 1, it can be assumed  that $G$ has an abelian
subgroup $A$ of finite index.

We claim that the commutator subgroup of $H=P\cdot A$ is a bounded
$p$-group.  Clearly $S_*(FP)$ satisfies a group identity and
according to Lemma 3  $P'$ is a bounded $p$-group. The normal
abelian $p$-subgroup $P'\cap A$ has  finite exponent and according
to Lemma 6 the ideal ${ \mathfrak{I}}(P'\cap A)$ is  locally
nilpotent of bounded degree. The subgroup $P'\cap A$ of $P'$ has
finite index in $P$ and
$$
{ \mathfrak{I}}(P')/ \mathfrak{I}(P'\cap A)\cong {
\mathfrak{I}}(P'/(P'\cap A)).
$$
Therefore  ${ \mathfrak{I}}(P')$ is a locally nilpotent ideal of
bounded degree $p^t$ for some $t$. Clearly  $FG/{
\mathfrak{I}}(P')\cong F[G/P']$ and  put  $P'=\gp{1}$. Since $A$
has a finite index in $H=P\cdot A$,  Lemma 3 ensures that $H'$ is
a $p$-group of bounded exponent and according to Lemma 1,  we can
put $H'=\gp{1}$  again.

The p-Sylow subgroup $P$ of $G$ is abelian and by Lemma 8 the
ideal  ${ \mathfrak{I}}(P)$ satisfies a non-matrix polynomial
identity, Moreover, by Lemma 9 the ideal  ${ \mathfrak{I}}(P)$
satisfies polynomial identities  of the following forms:\quad
$[x,y]^{p^l}$ and\quad $([x,y]z)^{p^l}$.

Let $h\in G$ and $a\in P$. Clearly  $(a-1)h, h^{-1}(a^{-1}-1)\in {
\mathfrak{I}}(P)$ and
$$
[(a-1)h,
h^{-1}(a^{-1}-1)]^{p^l}=(a^h)^{p^l}+(a^h)^{-p^l}-a^{p^l}-a^{-p^l}=0
$$
which implies  that either \quad $(h,a)^{p^l}=1$ or \quad
$h^{-1}a^{p^l}h=a^{-p^l}$.

Put $z=a^{p^l}$. From\quad   $h^{-1}a^{p^l}h=a^{-p^l}$\quad  it
follows that \quad $h^{-1}zh=z^{-1}$ and\qquad
$([z-1,(z^{-1}-1)h])^{p^l}=0$. Clearly \quad
$[z-1,(z^{-1}-1)h]=-z^{-2}(z+1)(z-1)^2h$\quad  so
\[
\begin{split}
0&=([z-1,(z^{-1}-1)h])^{p^l}\\
&=-\big((z+1)(z-1)^2(z^{-1}+1)(z^{-1}-1)^2h^2\big)^{\frac{p^l-1}{2}}\big(z^{-2}(z+1)(z-1)^2 h\big)\\
&=-z^{\frac{-3p^l-1}{2}}\cdot (z+1)^{p^l}\cdot (z-1)^{2p^l}\cdot
h^{p^l}.
\end{split}
\]
Since $char(K)>2$, the element $z+1$ is a unit and
$(z-1)^{2p^l}=(a-1)^{2p^{2l}}=0$ and the order of  $a$ at most
$2p^{2l}$. Therefore $(h,a)^{p^{2l+1}}=1$ for all $h\in G$, $a\in
P$ and $2l+1$ depends on only  the group identity. Since $(G,P)$
is a $p$-group of bounded exponent, we can again make a reduction,
so we can assumed that  $(G,P)=1$ and $P$ is central.

Let $P$ be a central subgroup of unbounded exponent and
$h_1,h_2\in G$. Obviously
\[
\begin{split}
((h_1,h_2)^{p^l}-1)(a-1)^{p^{3l}}&=((h_1,h_2)-1)^{p^l}(a-1)^{p^{3l}}\\
 &=([h_1^{-1}(a-1),h_2^{-1}(a-1)]h_1h_2(a-1))^{p^l}=0
\end{split}
\]
for $a\in P$. Since there are infinitely many element of the form
$a^{p^{3l}}$ we conclude that $(h_1,h_2)^{p^l}=1$ and the proof is
complete.
\end{proof}

%%%%%%%%%%%%%%%%%%%%%%%%%%%%%%%%%%%%%
\nocite{*}
\bibliographystyle{abbrv}
\bibliography{viktor0717}

\bigskip

\centerline{\em Received May 15, 2005.}  \centerline{\mbox{}}
\end{document}